\documentstyle[11pt]{article}
\def\Bbb R{{\rm \bf R}}
\def\proclaim#1{\vskip2mm{\bf #1}\em}
\def\endproclaim{\em \vskip2mm}
\def\tag#1{\eqno(#1)}
\def\gathered{\begin{array}{c}}
\def\endgathered{\end{array}}
\def\text{\mbox}

\begin{document}

\title {An approach based on the wave equation in the time domain for active shielding of an unwanted wave with a fixed frequency}
\author{Masaru IKEHATA\footnote{
Laboratory of Mathematics,
Graduate School of Advanced Science and Engineering,
Hiroshima University, Higashhiroshima 739-8527, JAPAN}
\footnote{Emeritus Professor at Gunma University}
}
%\date{}
\maketitle

\begin{abstract}
An approach for shielding an unwanted wave with a fixed frequency by generating a suitably controlled nontrivial wave with the same frequency
is suggested. Unlike the well known surface potential approach, the source of the controlled wave is given by solving the Cauchy problem for the wave equation in the finite time domain.

%{\bf Shield.tex}

\noindent
AMS:  35J05, 35L05, 35C15, 31B10, 35R30

\end{abstract}

%\tableofcontents

\section{Introduction}

Assume that we are hearing an unwanted sound wave with a known fixed frequency caused by a known source.
The support of the source is contained in the closure of a bounded open subset of $\Bbb R^3$ which we denote by $D$.  
We are standing outside a known bounded domain $\Omega$ that contains $\overline D$. 
The problem considered in this paper is: can one shield the wave outside domain $\Omega$ {\it completely}, by adding another wave
generated by a controlled source which is supported in $\overline\Omega$ and different from that of the unwanted wave?

This is a typical problem of an active shielding of a given unwanted wave.
This note is concerend with the methodology for
the {\it virtual sound barriers} \cite{Q} which is a special form of active noise control systems and avoids 
to block air, light, and access for shielding.
We suggest one natural approach which is based on the wave equation in the finite time domain.
To our best knowledge, our approach is new and, in particular, not listed in the book \cite{Q}, see pages 26-33 therein.

Now let us formulate the problem more precisely.
We assume that the unwanted wave satisfies 
the inhomogeneous Helmholtz equation
$$\begin{array}{ll}
\displaystyle
(\Delta+k^2)w+F(x)=0, & x\in\Bbb R^3
\end{array}
$$
and the outgoing radiation condition
$$\begin{array}{lll}
\displaystyle
\lim_{r\rightarrow\infty}r\left(\frac{\partial}{\partial r}-ik\right)w(r\omega)=0, & 
\displaystyle
r=\vert x\vert, 
& 
\displaystyle
\omega=\frac{x}{\vert x\vert},
\end{array}
$$
where $k$ is a positive number, $F\in L^2(\Bbb R^3)$ with $\text{supp}\,F\subset\overline D$
and the limit is uniform with respect to $\omega$.  It is known that, by virtue of the radiation condition,
$w$ has the expression
$$\begin{array}{ll}
\displaystyle
w(x)=\frac{1}{4\pi}\int_D\frac{e^{ik\vert x-y\vert}}{\vert x-y\vert}\,F(y)dy, & x\in\Bbb R^3.
\end{array}
\tag {1.1}
$$
See \cite{CK}.

We consider the following problem.

{\bf\noindent Problem.}  Find another source term $G\in L^2(\Bbb R^3)$ whose support is contained in $\overline\Omega$ and satisfies
$$\displaystyle
\text{supp}\,G\cap\overline D=\emptyset
\tag {1.2}
$$ 
in such a way that the solution of
the inhomogeneous Helmholtz equation
$$\begin{array}{ll}
\displaystyle
(\Delta+k^2)\tilde{w}+G(x)=0, & x\in\Bbb R^3
\end{array}
$$
with the outgoing radiation condition
$$\begin{array}{lll}
\displaystyle
\lim_{r\rightarrow\infty}r\left(\frac{\partial}{\partial r}-ik\right)\tilde{w}(r\omega)=0, & 
\displaystyle
r=\vert x\vert, 
& 
\displaystyle
\omega=\frac{x}{\vert x\vert},
\end{array}
$$
satisfies 
$$\begin{array}{ll}
\displaystyle
w(x)+\tilde{w}(x)=0,
&
\displaystyle
x\in\Bbb R^3\setminus\overline\Omega.
\end{array}
\tag {1.3}
$$
The $G$ may depend on $k$, $D$, $F$ and $\Omega$ and is called an active source.
The $\tilde{w}$ which is called the secondary wave, has the expression
$$\begin{array}{ll}
\displaystyle
\tilde{w}(x)=\frac{1}{4\pi}\int_{\Omega\setminus D}\frac{e^{ik\vert x-y\vert}}{\vert x-y\vert}\,G(y)dy, & x\in\Bbb R^3.
\end{array}
\tag {1.4}
$$

Since we require the condition (1.2), one can not choose the {\it trivial solution} $G=-F$.
And in practice, condition (1.2) is natural, say, consider the case when the unwanted wave is radiated from
something {\it volumemetric real body}.

In this note, to clearly indicate the idea we consider the simplest case when $D=B_{\epsilon}$, where $B_{\epsilon}$ is an open
ball with radius $\epsilon$ and centered at the origin of the Cartesian coordinates.

\section{A solution}

The conclusion is: one can construct the active source as a function of the source of the unwanted wave
by using a solution of the wave equation in the finite time domain.

Let us describe the solution $G$ to the problem in the case $D=B_{\epsilon}$ step by step.

(1)  Let $T>2\epsilon$.
Solve the Cauchy problem for the wave equation in the whole space
$$\left\{
\begin{array}{ll}
\displaystyle
(\partial_t^2-\Delta)u=0, & x\in\Bbb R^3,\,0<t<T,\\
\\
\displaystyle
u(x,0)=0, & x\in\Bbb R^3,\\
\\
\displaystyle
\partial_tu(x,0)=F(x), & x\in\Bbb R^3.
\end{array}
\right.
\tag {2.1}
$$

(2)  Define
$$\begin{array}{ll}
\displaystyle
G(x)=-e^{-ikT}(\partial_tu(x,T)+iku(x,T)), & x\in\Bbb R^3.
\end{array}
\tag {2.2}
$$

Since $\text{supp}\,F\subset\overline B_{\epsilon}$, by Huygens's principle, we have 
$$\displaystyle
\text{supp}\,G\subset\overline{B_{T+\epsilon}}
$$
and
$$\displaystyle
\text{supp}\,G\subset\Bbb R^3\setminus B_{T-\epsilon}.
\tag {2.3}
$$

Thus the support of $G$ is contained in the shell domain $\overline B_{T+\epsilon}\setminus B_{T-\epsilon}$ and satisfies (1.2) provided $T>2\epsilon$.
See Figure 1 for an illustration of a situation.

Note that (2.3) is a consequence of the fact that the dimension of $\Bbb R^3$ is an odd number.
See, e.g., \cite{ES}.  

Thus, formula (1.4) becomes
$$\begin{array}{ll}
\displaystyle
\tilde{w}(x)=\frac{1}{4\pi}\int_{B_{T+\epsilon}\setminus B_{T-\epsilon}}\frac{e^{ik\vert x-y\vert}}{\vert x-y\vert}\,G(y)dy, & x\in\Bbb R^3.
\end{array}
\tag {2.4}
$$

Now we are ready to state the main part of this note.

\proclaim{\noindent Theorem 2.1.}  Let $D=B_{\epsilon}$ and choose $\Omega=B_{T+\epsilon}$.
Let $T>2\epsilon$.  Then,  the $G$ given by (2.2) satisfies (1.3).

\endproclaim

{\it\noindent Proof.}
Define
$$\begin{array}{ll}
\displaystyle
z(x)=\int_0^T e^{-ikt}u(x,t)\,dt, & x\in\Bbb R^3.
\end{array}
$$
It follows from (2.1) that the $z$ satisfies
$$\begin{array}{ll}
\displaystyle
\Delta z+k^2 z+F(x)+G(x)=0, & x\in\Bbb R^3.
\end{array}
$$
And we have
$$\displaystyle
\text{supp}\,z\subset\overline{B_{T+\epsilon}}.
\tag {2.5}
$$
These imply that $z$ has the expression
$$\displaystyle
z=w+\tilde{w}.
$$
Now (2.5) yields the desired conclusion (1.3).

\noindent
$\Box$

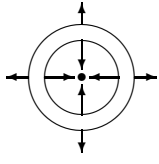
\begin{figure}
\setlength\unitlength{0.1truecm}
\begin{picture}(10,7)(-5,1)
\put(65,10){\circle*{0.9}{}}
\put (65,10){\circle {20}}
\put (65,10){\circle {10}}
\put(72,10){\vector(1,0){3}}
\put(70,10){\vector(-1,0){4}}
  \put(65,17){\vector(0,1){3}}
  \put(65,15){\vector(0,-1){4}}
  \put(65,3){\vector(0,-1){3}}
  \put(65,5){\vector(0,1){4}}
  \put(58,10){\vector(-1,0){3}}
  \put(60,10){\vector(1,0){4}}
\end{picture}
\caption{\label{fig:1}The shell domain surrounding $B_{\epsilon}$.}
\end{figure}

Note that $G$ depends also on $T$. Of course,  $T$ sholud not be so large since $\Omega=B_{T+\epsilon}$ and $\Omega$ becoms large
if $T$ is large.

By Theorem 2.1, we know that one can completely
shield the wave field with known source $F$ by adding another wave field.
The added wave field has a source whose support is located outside of $\text{supp}\,F$ and thus is nontrivial.
The proof employs a solution of the Cauchy problem for the wave equation over a {\it finite time interval}.

As far as the author knows, no one has pointed out this fact.  It should be pointedt out that
the construction of $G$ which is given by (2.2) was inspired by reconsidering the {\it time domain enclosure method} \cite{IE0, IW2, IW3, IW4, IW5}.

\section{A consequence on the far field pattern}

From (1.1) and (1.4) we see that functions $w$ and $\tilde{w}$
have the expressions as $r\rightarrow\infty$ for all $\omega\in S^2$ uniformly
$$\displaystyle
w(r\omega)\sim\frac{1}{4\pi}\frac{e^{ikr}}{r}\int_De^{-iky\cdot\omega}F(y)dy
$$
and
$$\displaystyle
\tilde{w}(r\omega)\sim\frac{1}{4\pi}\frac{e^{ikr}}{r}\int_{B_{T+\epsilon}\setminus B_{T-\epsilon}}e^{-iky\cdot\omega}G(y)dy.
$$
Therefore Theorem 2.1 ensures that
$$\displaystyle
\int_De^{-iky\cdot\omega}F(y)dy+\int_{B_{T+\epsilon}\setminus B_{T-\epsilon}}e^{-iky\cdot\omega}G(y)dy=0
$$
under the choice of $G$ given by (2.2).  This means that the far field pattern (see \cite{CK}) of $w+\tilde{w}$ vanishes.

In particular, the far field pattern of $w$ coincides with that of $-\tilde{w}$.
This gives an example for the nonuniqueness for the inverse source problem: one can not uniquely determine the source term 
for the Helmholtz equation at a fixed frequency from the far-field pattern.

This fact itself is well known, however, this example tells us more than a non uniqueness in inverse source problem.
One can {\it hide} the field radiated from a known source by generating an another field radiated by a set of suitable sources distributed
around the known source.

\section{Comparison with an approach based on surface potentials}

In this section we present another approach in \cite{LRT} which is based on surface potentials.
In more gerneral settings see also \cite{U}.  And for references in early studies see also introduction in \cite{NU}.  

The idea is simple.

Set
$$\begin{array}{ll}
\displaystyle
K(y-x)=\frac{1}{4\pi}\frac{e^{ik\vert y-x\vert}}{\vert y-x\vert}, & x\not=y.
\end{array}
$$

Let $x\in\Bbb R^3\setminus\overline\Omega$.

From the governing equation of $w$ in $\Bbb R^3\setminus\overline\Omega$
and the radiation condition, we have the expression
$$\begin{array}{ll}
\displaystyle
w(x)=\int_{\partial\Omega}\left(w\frac{\partial K_x}{\partial\nu}-
K_x\frac{\partial}{\partial\nu}w\right)\,dS(y), & x\in\Bbb R^3\setminus\overline\Omega,
\end{array}
\tag {4.1}
$$
where $K_x(y)=K(y-x)$.
Note that the support of $F$ is contained in $\Omega$.  This is a well known formula in scattering theory, see e.g., \cite{CK}.

Rewrite (4.1) as
$$\begin{array}{ll}
\displaystyle
w(x)-\int_{\partial\Omega}\left(w\frac{\partial K_x}{\partial\nu}-
K_x\frac{\partial}{\partial\nu}w\right)\,dS(y)=0, & x\in\Bbb R^3\setminus\overline\Omega.
\end{array}
\tag {4.2}
$$
In \cite{LRT} they consider that this is a cancelation formula of $w$ outside $\Omega$.  They define instead of $\tilde{w}$
$$\begin{array}{ll}
\displaystyle
w'(x)=-\int_{\partial\Omega}\left(w\frac{\partial K_x}{\partial\nu}-
K_x\frac{\partial}{\partial\nu}w\right)\,dS(y), & x\in\Bbb R^3.
\end{array}
\tag {4.3}
$$
Then, (4.2) means that
$$\begin{array}{ll}
\displaystyle
w(x)+w'(x)=0, & x\in\Bbb R^3\setminus\overline\Omega.
\end{array}
$$

The remarkable point of (4.3) is: the $w'$ depends on only the Cauchy data of $w$ on $\partial\Omega$.
It does not require any detailed knowledge of the source of $w$ in $\Omega$.

Note also that $\Omega$ can be an arbitrary bounded domain.

In \cite{Q}, on page 28 it is clarified that the theoretical base of the virtual sound barriers
is such type of formulae which is called the Kirchhoff-Helmholtz equation and a quantified version of
Huygens's principle.
Mathematically, it is an application of integration by parts combined with
the property of the fundamental solution of the Helmholtz equation.

Following \cite{LRT}, one can rewrite (4.3) more.

Choose a function $\Psi\in C_0^2(\Bbb R^3)$ in such a way that
$$\begin{array}{lll}
\displaystyle
\Psi=w, & 
\displaystyle
\frac{\partial \Psi}{\partial\nu}=\frac{\partial w}{\partial\nu}, 
& x\in\partial\Omega.
\end{array}
\tag {4.4}
$$
We have, for all $x\in\Bbb R^3\setminus\overline\Omega$
$$\begin{array}{l}
\displaystyle
\,\,\,\,\,\,
\int_{\partial\Omega}\left(w\frac{\partial K_x}{\partial\nu}-
K_x\frac{\partial}{\partial\nu}w\right)\,dS(y)\\
\\
\displaystyle
=\int_{\partial\Omega}\left(\Psi\frac{\partial K_x}{\partial\nu}-
K_x\frac{\partial}{\partial\nu}\Psi\right)\,dS(y)
\\
\\
\displaystyle
=\int_{\Bbb R^3\setminus\Omega}K(y-x)(\Delta+k^2)\Psi(y)\,dy-\Psi(x)
\end{array}
$$
and
$$\begin{array}{ll}
\displaystyle
\Psi(x)=\int_{\Bbb R^3}K(y-x)(\Delta+k^2)\Psi(y)\,dy, & x\in\Bbb R^3.
\end{array}
$$

Thus, we obtain, for all $x\in\Bbb R^3\setminus\overline\Omega$
$$\begin{array}{ll}
\displaystyle
\,\,\,\,\,\,
-\int_{\partial\Omega}\left(w\frac{\partial K_x}{\partial\nu}-
K_x\frac{\partial}{\partial\nu}w\right)\,dS(y)
\\
\\
\displaystyle
=\int_{\Omega}K(y-x)(\Delta+k^2)\Psi(y)dy.
\end{array}
$$
Therefore, we obtain
$$\begin{array}{ll}
\displaystyle
w'(x)=\int_{\Omega}K(y-x)(\Delta+k^2)\Psi(y)dy, & x\in\Bbb R^3\setminus\overline\Omega.
\end{array}
$$

Next let $x\in\Omega$.  We have
$$\begin{array}{l}
\displaystyle
\,\,\,\,\,\,
-\int_{\partial\Omega}\left(w\frac{\partial K_x}{\partial\nu}-
K_x\frac{\partial}{\partial\nu}w\right)\,dS(y)\\
\\
\displaystyle
=-\int_{\partial\Omega}\left(\Psi\frac{\partial K_x}{\partial\nu}-
K_x\frac{\partial}{\partial\nu}\Psi\right)\,dS(y)
\\
\\
\displaystyle
=\int_{\Bbb R^3\setminus\Omega}K(y-x)(\Delta+k^2)\Psi(y)\,dy
\end{array}
$$
Therefore, we obtain
$$\begin{array}{ll}
\displaystyle
w'(x)=\int_{\Bbb R^3\setminus\Omega}K(y-x)(\Delta+k^2)\Psi(y)dy, & x\in\Omega.
\end{array}
$$

Summing up, we have obtained the following formula.

\proclaim{\noindent Proposition 4.1.}
Let $\Psi\in C^2_0(\Bbb R^3)$ satisfy (4.4).
Then, the $w'$ given by (4.3) has the the expression
$$\displaystyle
w'(x)=
\left\{
\begin{array}{ll}
\displaystyle
\int_{\Bbb R^3\setminus\Omega}K(y-x)(\Delta+k^2)\Psi(y)dy, & x\in\Omega,\\
\\
\displaystyle
\int_{\Omega}K(y-x)(\Delta+k^2)\Psi(y)dy, & x\in\Bbb R^3\setminus\overline\Omega.
\end{array}
\right.
\tag {4.5}
$$

\endproclaim

Note that the support of $\Psi$ can be an arbitrary small neighbourhood of $\partial\Omega$.
One can rewrite (4.5) compactly
$$\displaystyle
w'=\chi_{\Omega}K*\left\{\chi_{\Bbb R^3\setminus\overline\Omega}(\Delta+k^2)\Psi\right\}
+\chi_{\Bbb R^3\setminus\overline\Omega}K*\left\{\chi_{\Omega}(\Delta+k^2)\Psi\right\}.
$$
We see that $w'\vert_{\Omega}\in H^2(\Omega)$, $w'\vert_{B_R\setminus\overline\Omega}\in H^2(B_R\setminus\overline\Omega)$ for all $R>>1$
and $w'\in L^2(\Bbb R^3)$.  

However, $w'$ itself does not belong to $H^1$ in a neighbourhood of $\partial\Omega$.
This is a consequence of the jump relation of the double layer potential and the radiation condition.
Thus $w'$ can not be realized as a wave field having a compact source in $L^2(\Bbb R^3)$.

In contrast to this, our construction of $\tilde{w}$ automatically ensures $\tilde{w}\in H^2_{\text{loc}}(\Bbb R^3)$.
This is an advantage of making use of the full knowledge of the source of the unwanted wave $w$.

Le us make a comparison.

\noindent
{\bf The surface potential method.}

$\bullet$  The selection of $\Omega$ is arbitrary as long as the condition $\text{suupp}\,F\subset\Omega$
is satisfied.   For the construction of $w'$ we need only the Cauchy data of $w$ on $\partial\Omega$.

$\bullet$  $w'$ has a {\it singularity across} $\partial\Omega$ and its source does not belong to $L^2(\Bbb R^3)$.

\noindent
{\bf Our method.}

$\bullet$  The selection of $\Omega$ depends on an upper bound of the size of $\text{supp}\,F$.  
The full knowledge of the source of the unwanted wave is required.

$\bullet$  $\tilde{w}$ is locally $H^2$-regular in the whole space and its source belongs to $L^2(\Bbb R^3)$.  
It may be possible to approximate the source as a superposition of finitely many {\it monopoles} only.
See (2.4).

\section{Comparison with a naive approach in a special case}

Consider the case $F(x)=(\epsilon-\vert x\vert)\chi_{B_{\epsilon}}(x)$.  The $F$ belongs to $H^1(\Bbb R^3)$.
The field generated by the source $F$ is given by
$$\displaystyle
w(x)=\frac{1}{4\pi}\int_{B_{\epsilon}}\frac{e^{ik\vert x-y\vert}}{\vert x-y\vert}\,(\epsilon-\vert y\vert)\,dy.
$$

Here we recall
\proclaim{\noindent Lemma 5.1.}
Let $B$ be the open ball with radius $\eta$ and centered at the origin. 
We have, for all $x\in\Bbb R^3\setminus\overline B$
$$\displaystyle
\frac{1}{4\pi}\int_B\frac{e^{ik\vert x-y\vert}}{\vert x-y\vert}\,dy
=-i\frac{\varphi(-ik\eta)}{k^3}\frac{e^{ik\vert x\vert}}{\vert x\vert}
\tag {5.1}
$$
and
$$\displaystyle
\frac{1}{4\pi}\int_B\frac{e^{ik\vert x-y\vert}}{\vert x-y\vert}\,(\eta-\vert y\vert)\,dy
=\frac{1}{k^4}\frac{e^{ik\vert x\vert}}{\vert x\vert}
P(-ik\eta),
\tag {5.2}
$$
where
$$\left\{
\begin{array}{l}
\displaystyle
\varphi(s)=s\cosh s-\sinh s,
\\
\\
\displaystyle
P(s)=-2\cosh s+s\sinh s+2.
\end{array}
\right.
$$
\endproclaim

The equation (5.1) is a consequence of the mean value theorem for the Helmholtz equation. 
For the proof see \cite{CH} and that of (5.2) see Appendix in \cite{IW4}.

Applying (5.2) to $w(x)$ for $x\in\Bbb R^3\setminus\overline {B_{\epsilon}}$, we obtain
$$\displaystyle
w(x)=\frac{1}{k^4}\frac{e^{ik\vert x\vert}}{\vert x\vert}
P(-ik\epsilon).
\tag {5.3}
$$

Consider also the secondary field $w'(x)$ generated by the source $G(x)=\chi_{B_{R_2}\setminus B_{R_1}}(x)$ with $R_2>R_1>\epsilon$:
$$\displaystyle
w'(x)=\frac{1}{4\pi}\int_{B_{R_2}\setminus B_{R_1}}\frac{e^{ik\vert x-y\vert}}{\vert x-y\vert}\,dy.
$$

From (5.1) we have, for all $x\in\Bbb R^3\setminus\overline{B_{R_2}}$
$$\displaystyle
w'(x)=-\frac{i}{k^3}\frac{1}{4\pi}\frac{e^{ik\vert x\vert}}{\vert x\vert}(\varphi(-ik R_2)-\varphi(-ik R_1)).
\tag {5.4}
$$

Thus a combination of equations (5.3) and (5.4) yields that: $w(x)+w'(x)=0$ for all $x\in\Bbb R^3\setminus B_{R_2}$ if and only if
$$\displaystyle
\frac{1}{k^4}
P(-ik\epsilon)-\frac{i}{k^3}(\varphi(-ik R_2)-\varphi(-ik R_1))=0,
$$
that is
$$\displaystyle
P(-ik\epsilon)-ik(\varphi(-ik R_2)-\varphi(-ik R_1))=0.
\tag {5.5}
$$

So, the problem is to: given $\epsilon$ and $k$ find $R_2>R_1>\epsilon$ such that (5.5) is valid.

We have
$$\begin{array}{l}
\displaystyle
\,\,\,\,\,\,
P(-ik\epsilon)\\
\\
\displaystyle
=-2\cosh (-ik\epsilon)-ik\epsilon\sinh (-ik\epsilon)+2
\\
\\
\displaystyle
=-2\cos k\epsilon-k\epsilon \sin k\epsilon+2,
\end{array}
$$
$$\begin{array}{ll}
\displaystyle
\,\,\,\,\,\,
\varphi(-ikR_j)
\\
\\
\displaystyle
=-ikR_j\cosh (-ikR_j)-\sinh(-ikR_j)\\
\\
\displaystyle
=-ikR_j\cos kR_j+i\sin kR_j
\\
\\
\displaystyle
=-i(kR_j\cos kR_j-\sin kR_j).
\end{array}
$$

Thus, (5.5) is equivalent to the equation
$$\displaystyle
-2\cos k\epsilon-k\epsilon \sin k\epsilon+2
=k(kR_2\cos kR_2-\sin kR_2-kR_1\cos kR_1+\sin kR_1).
\tag {5.6}
$$

Define
$$\displaystyle
Q(\xi)=\xi\cos\xi-\sin\xi.
$$
Then, (5.6) becomes
$$\displaystyle
Q(kR_2)=Q(kR_1)+\frac{1}{k}\left\{2(1-\cos k\epsilon)-k\epsilon\sin k\epsilon\right\}.
\tag {5.7}
$$

From the behaviour of $Q(\xi)$ we see that, given $k>0$, $\epsilon>0$ and $R_1>\epsilon$ there exist {\it infinitely many} $R_2>R_1$ such that
(5.7) is satisfied.  Thus, the shielding is possible and the radius $R_2$ is characterized as a soultion of the equation (5.7) with $R_2>R_1$.

However, if the source term has the form 
$$\displaystyle
F(x)=\rho(x)\chi_{B_{\epsilon}}(x),
$$
where $\rho(x)$ is a function having a general form, the construction of $w'$ having the form
$$\begin{array}{ll}
\displaystyle
w'(x)=\frac{1}{4\pi}\int_{B_{R_2}\setminus B_{R_1}}\frac{e^{ik\vert x-y\vert}}{\vert x-y\vert}\,g(y)\,dy,
&
R_2>R_1>\epsilon,
\end{array}
$$
shall be difficult without using our natural approach based on Huygens's principle.

Our method gives a simple solution
that $R_2=T+\epsilon$, $R_1=T-\epsilon$ with $T>2\epsilon$ and $g(y)$ is given by the right-hand side on (2.2).
This means that we made use of the Cauchy problem for wave equation (2.1) as a {\it calculater} of the desired source.

\section{Testing the approach}

It would be interesting to do numerical testing of our method based on formula (2.4).
The method consists of only two steps.

$\quad$

{\bf Step 1.}  Give $F$ and compute $G$ given by (2.2).

$\quad$

{\bf Step 2.}  Generate the secondary wave $\tilde{w}$ with the source $G$ computed in Step 1.

$\quad$

In Step 1, the computation needs the solution of the Cauchy problem for the wave equation (2.1)
together with its time derivative at $t=T$, however, not for all $t\in\,]0,\,T[$.  
We have the Kirchhoff formula, e.g., \cite{ES} for the solution of (2.1) which yields the exact value
of the solution together with the time derivative at $t=T$ without solving (2.1) numerically.

Numerically, Step 2 means that: compute $\tilde{w}$
via formula (2.4).

To check our method numerically, we have to compute also the unwanted wave $w$ via formula (1.1).
Then, compute the total wave $w+\tilde{w}$ and observe its behaviour for $x\in B$,
where $B$ is a ball centered at the origin with a large radius.

$$\quad$$

\centerline{{\bf Acknowledgment}}

The author was partially supported by Grant-in-Aid for
Scientific Research (C)(No. 17K05331) of Japan  Society for
the Promotion of Science.

$$\quad$$

\vskip1cm
\noindent
e-mail address

ikehata@hiroshima-u.ac.jp

\end{document}